\documentclass[reqno]{amsart}
\usepackage{cite,graphicx,calrsfs}
\usepackage{xcolor}

\usepackage{caption}
\theoremstyle{definition}
\newtheorem{theorem}{Theorem}

\theoremstyle{remark}
\newtheorem*{remark}{Remark}

\makeatletter 
\def\thmhead@plain#1#2#3{%
  \thmname{#1}\thmnumber{\@ifnotempty{#1}{ }\@upn{#2}}%
  \thmnote{{\the\thm@notefont#3}}}
\let\thmhead\thmhead@plain
\makeatother

\begin{document} 
\author{Dmitriy Dmitrishin}
\address{Dmitriy Dmitrishin: Odessa National Polytechnic University, 1 Shevchenko Ave., Odesa 65044, Ukraine}
\email{dmitrishin@op.edu.ua}

\author{Alexander Stokolos}
\address{Alexander Stokolos: Georgia Southern University, Statesboro GA, 30460, USA }
\email{astokolos@georgiasouthern.edu}
\date{\today}

\title[Some properties of quadrinomials]{Some properties of the quadrinomials
$p(z)=1+\kappa(z+z^{N-1})+z^N$ and $q(z)=1+\kappa(z-z^{N-1})-z^N$}


\begin{abstract}\vspace*{1cm}
We show that all the zeros of the quadrinomial $p(z)=1+\kappa(z+z^{N-1})+z^N$ lie on the 
unit circle if and only if the inequalities 
\[
-1\le\kappa\le
\left\{\begin{array}{ll}
     1 & \text{ if $N$ is {\it even}}, \\
     N/(N-2) & \text{ if $N$ is {\it odd}}
\end{array}\right.
\]
hold. For the quadrinomial $q(z)=1+\kappa(z-z^{N-1})-z^N$, the corresponding
inequalities are 
\[
-N/(N-2)\le\kappa\le
\left\{\begin{array}{ll}
     1 & \text{ if $N$ is {\it odd}}, \\
     N/(N-2) & \text{ if $N$ is {\it even}}. 
\end{array}\right.
\]
In the cases of limiting values of the parameter $\kappa$, we provide factorization 
formulas for the corresponding quadrinomials. For example, when $N$ is {\it odd} and 
$\kappa=N/(N-2)$, the following representation is valid:
\[
p(z)=(1+z)^3\prod_{j=1}^{(N-3)/2}[1+z^2-2z\gamma_j],
\]
where $\gamma_j=1-2\nu_j^2$ with $\{\nu_j\}_{j=1}^{(N-3)/2}$ being the collection of 
positive roots of the equation $U'_{N-2}(x)=0$; here
\[
U_j(x)=U_j(\cos t)=\frac{\sin(j+1)t}{\sin t}=2^j x^j+\ldots
\]
are Chebyshev polynomials of the second kind and $U'_j(x)$ are their derivatives. Similar
factorization formulas are also provided for $q(z)$. As an application of the obtained 
results, we give the factorization formulas for the derivative of the Fej\'er polynomial, 
as well as construct certain univalent polynomials related to the polynomials $p(z)$ and 
$q(z)$.
\end{abstract} 
\maketitle 

\vspace{-11cm}
\hspace*{6cm }{\it\small In memory of Vladimir Dmitrishin \vspace*{4cm}}\\
\vspace{5.5cm}
\section{Introduction}
\begin{figure}[ht]
{\footnotesize
\begin{minipage}{.4\textwidth}
\begin{center}
\includegraphics[width=\textwidth]{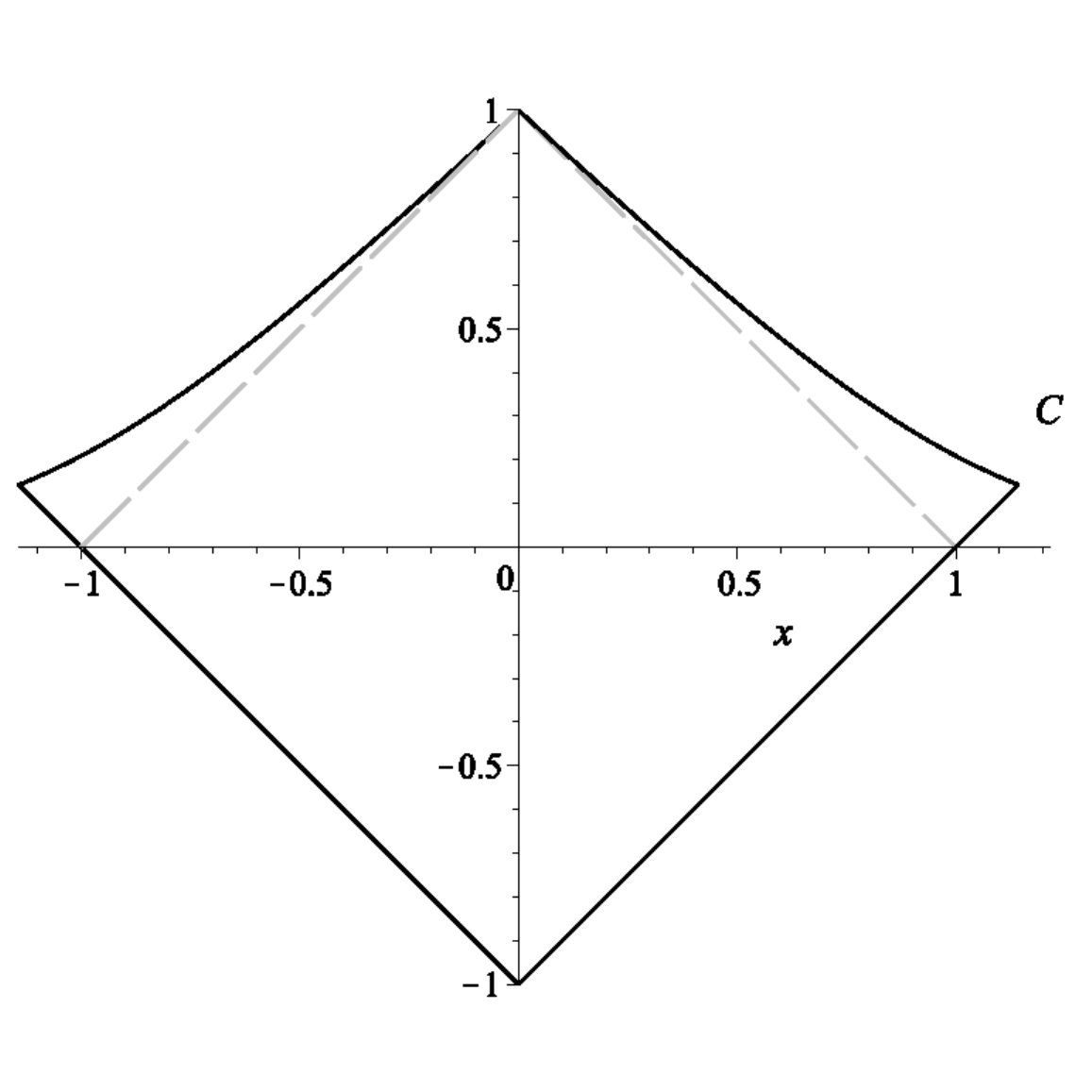}\\
\textit{i)}
\end{center}
\end{minipage}
\hspace{12pt}
\begin{minipage}{.4\textwidth}
\begin{center}
\includegraphics[width=\textwidth]{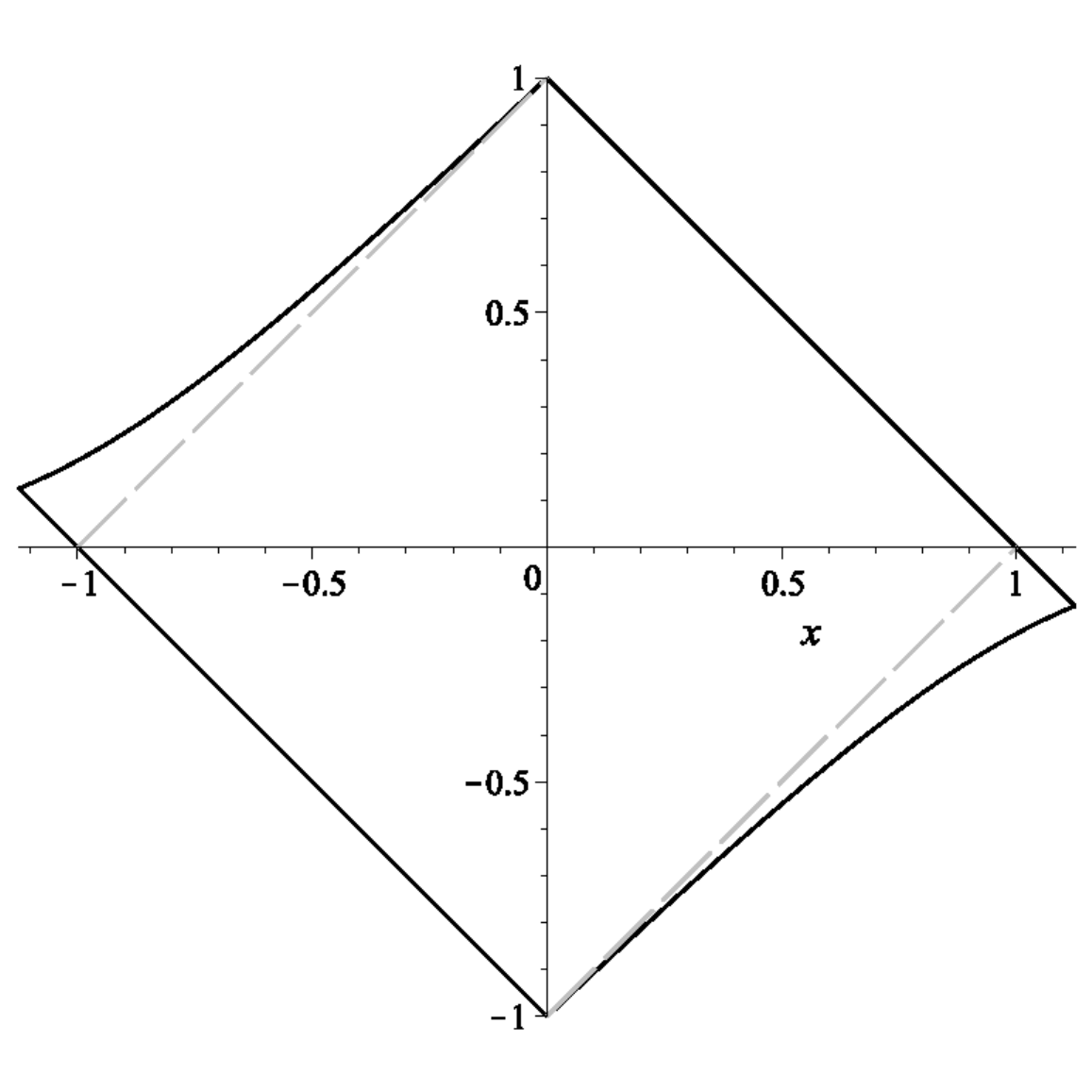}\\
\textit{ii)}
\end{center}
\end{minipage}
}
\caption{Stability domains of the trinomial $f(z)=z^n+az^{n-1}+b$ in the $(a,b)$-plane: 
{\it i)} corresponds to the case of even $n$, {\it ii)} is for odd $n$.}
\label{fig:domains}
\end{figure}

When working on {\it C.~Michel}'s problem of stretching the unit disk~\cite{CM71,DSS23}
and, correspondingly, on generalization of Theorem~\ref{th:intro-1} about representation
of typically real polynomials~\cite{MB89}, we came across the polynomial 
\begin{equation}\label{eq:curious-pln}
p(z)=1+\frac{N}{N-2}(z+z^{N-1})+z^N,
\end{equation}
which in some sense is an analogue of the polynomial $p(z)=1+z^N$. It turned out that
polynomial~\eqref{eq:curious-pln} has interesting properties, such as: it is factorized through 
the zeros of the derivative of a Chebyshev polynomial of the second kind with corresponding
order; it generates, by means of {\it T.~J.~Suffridge}'s transformation~\cite[p.~227]{TS72},
the univalent in the central unit disk polynomial
\[
F(z)=\frac{N}{N-1}\tilde{p}(z)-\frac{1}{N-1}z\tilde{p}'(z),
\]
where
\[
\tilde{p}(z)=\frac{z}{(1+z)^2}\left(1+\frac{N}{N-2}(z+z^{N+1})+z^N\right).
\]
All zeros of the polynomials~\eqref{eq:curious-pln} and $p(z)=1+z^N$ lie on the unit circle. 
Numerous studies are devoted to the polynomials with all the zeros on the unit circle, 
among which, e.g., \cite{TS76,EG91,WC95,KP08,DK11,PL02}. The following classical criterion
by {\it A.~Cohn} for the zeros of a polynomial belonging to the unit circle is well known. 
\begin{theorem}[{\cite[p.~121]{AC22}}]\label{th:intro-1}
All zeros of the polynomial 
\[
P(z)=\sum_{j=0}^N a_j z^{N-j},\quad a_0\ne0,
\]
lie on the unit circle if and only if the polynomial $P(z)$ is self-reciprocal---which means,
for a polynomial with real coeficients, that $P(z)=\pm z^N P(1/z)$---and all zeros of the
polynomial $P'(z)$ belong to the closed central unit disk 
$\overline{\mathbb{D}}=\{z:|z|\le1\}$.
\end{theorem}
In the general case, criteria for polynomial zeros belonging to the closed disk 
$\overline{\mathbb{D}}$ or the open disk $\mathbb{D}=\{z:|z|<1\}$ are quite cumbersome. In 
some special cases, however, these criteria become much simpler.
\begin{theorem}[\cite{SK94}]\label{th:intro-2}
All zeros of the trinomial $f(z)=z^n+az^{n-1}+b$ belong to $\mathbb{D}$ if and only if
the pair $(a,b)$ is an interior point of the finite domain the boundaries of which are defined
by the four curves:
\begin{enumerate}
    \item[I.] $b=-a-1$, \ $a\in\left[-\frac{n}{n-1},0\right]$, 
    \item[II.] $b=(-1)^n(a-1)$, \ $a\in\left[0,\frac{n}{n-1}\right]$, 
    \item[III.] $\left\{\begin{array}{l}
         \displaystyle a=-\frac{\sin nt}{\sin(n-1)t}, \\
         \displaystyle b=\frac{\sin t}{\sin(n-1)t}, 
    \end{array}\right.
    t\in\left[\frac{n-1}{n}\pi,\pi\right)$,
    \item[IV.] $\left\{\begin{array}{l}
         \displaystyle a=\frac{\sin nt}{\sin(n-1)t}, \\
         \displaystyle b=(-1)^n\frac{\sin t}{\sin(n-1)t},
    \end{array}\right.
    t\in\left[\frac{n-1}{n}\pi,\pi\right)$. 
\end{enumerate}
\end{theorem}
The domains defined in Theorem~\ref{th:intro-2} are called the {\it stability domains} of the
trinomial $f(z)=z^n+az^{n-1}+b$ in the coefficient plane $(a,b)$. These domains are 
represented in Figure~\ref{fig:domains}.

A point $C$ lies on line~II when $a=\frac{n}{(n-1)}$, or line~III as $t\to\pi$; that is,
$C\left(\frac{n}{n-1},\frac{1}{n-1}\right)$. Thus, one may use 
Theorems~\ref{th:intro-1} and \ref{th:intro-2} to obtain criteria for the zeros of some
quadrinomials to belong to the unit circle. 

In~\cite{TS72}, {\it T.~J.~Suffridge} considered the class $\mathcal{P}_n$ of polynomials 
$f(z)=z+\sum_{j=2}^n a_j z^j$ of degree $n$ satisfying the following condition: the equations
\[
\frac{f(ze^{i\alpha_k})-f(ze^{-i\alpha_k})}{z(e^{i\alpha_k}-e^{-i\alpha_k})}=0,
\quad k=1,\ldots,n,
\]
where $\alpha_k=k\pi/(n+1)$, do not have solutions in $\mathbb{D}$. It is not difficult to 
check that 
\[
\frac{f(ze^{i\alpha_k})-f(ze^{-i\alpha_k})}{z(e^{i\alpha_k}-e^{-i\alpha_k})}
=1+\sum_{j=2}^n a_j\frac{\sin j\alpha_k}{\sin\alpha_k}z^{j-1}.
\]
The transform 
\[
f^*(z)=\frac{n+1}{n}f(z)-\frac1n zf'(z)
\]
was introduced in~\cite[p.~227]{TS72}; clearly, 
$f^*(z)=\sum_{j=1}^n \left(1-\frac{j-1}{n}\right)a_j z^j$.
\begin{theorem}[{\cite[p.~228]{TS72}}]\label{th:Suffridge}
If $f(z)\in\mathcal{P}_n$ and $|a_n|=1$, then $f^*(z)$ is univalent in 
$\mathbb{D}$.
\end{theorem}
\section{Main results}
\begin{figure}[ht]
{\footnotesize
\begin{minipage}{.45\textwidth}
\begin{center}
\includegraphics[width=\textwidth]{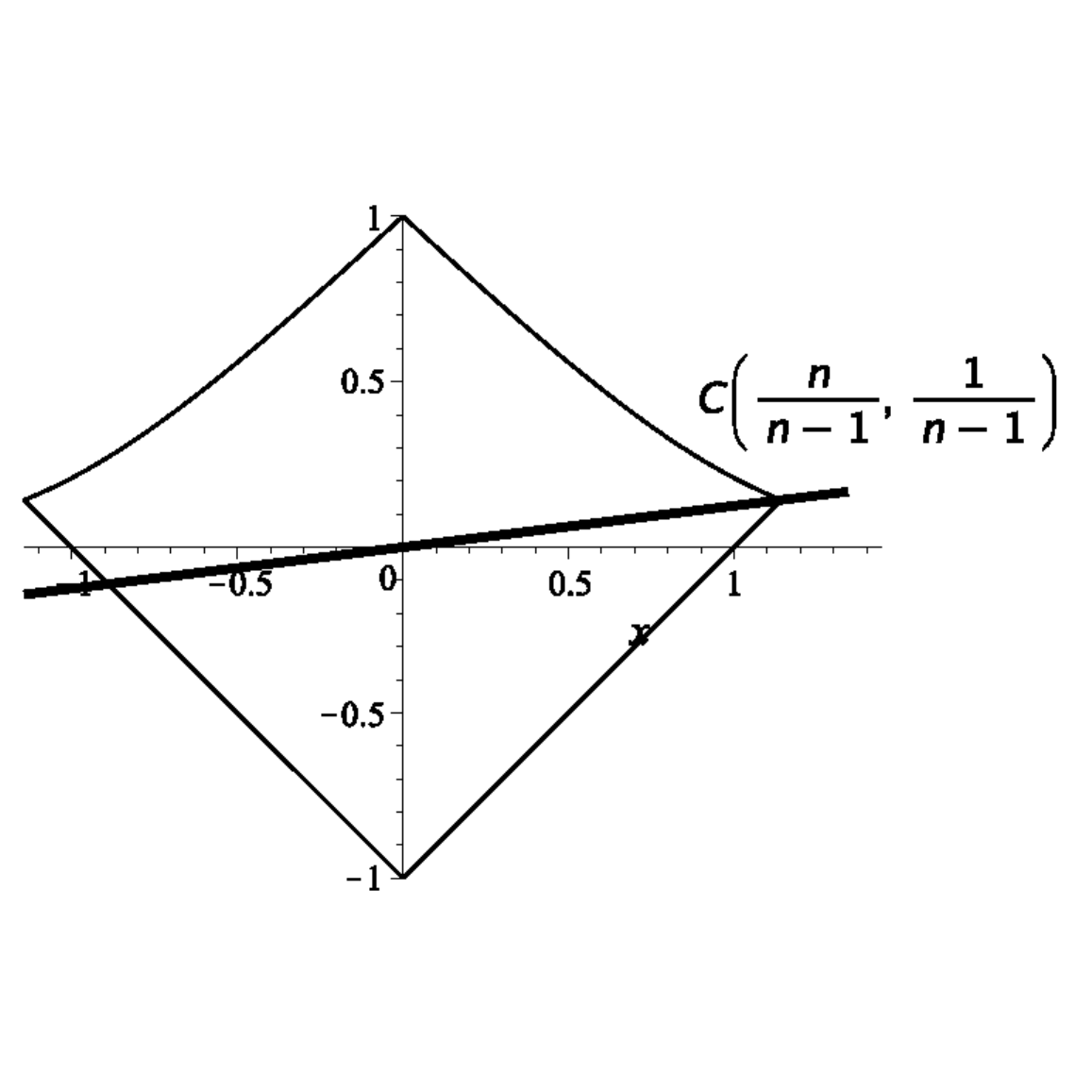}\\
\textit{i)}
\end{center}
\end{minipage}
\hspace{12pt}
\begin{minipage}{.35\textwidth}
\begin{center}
\includegraphics[width=\textwidth]{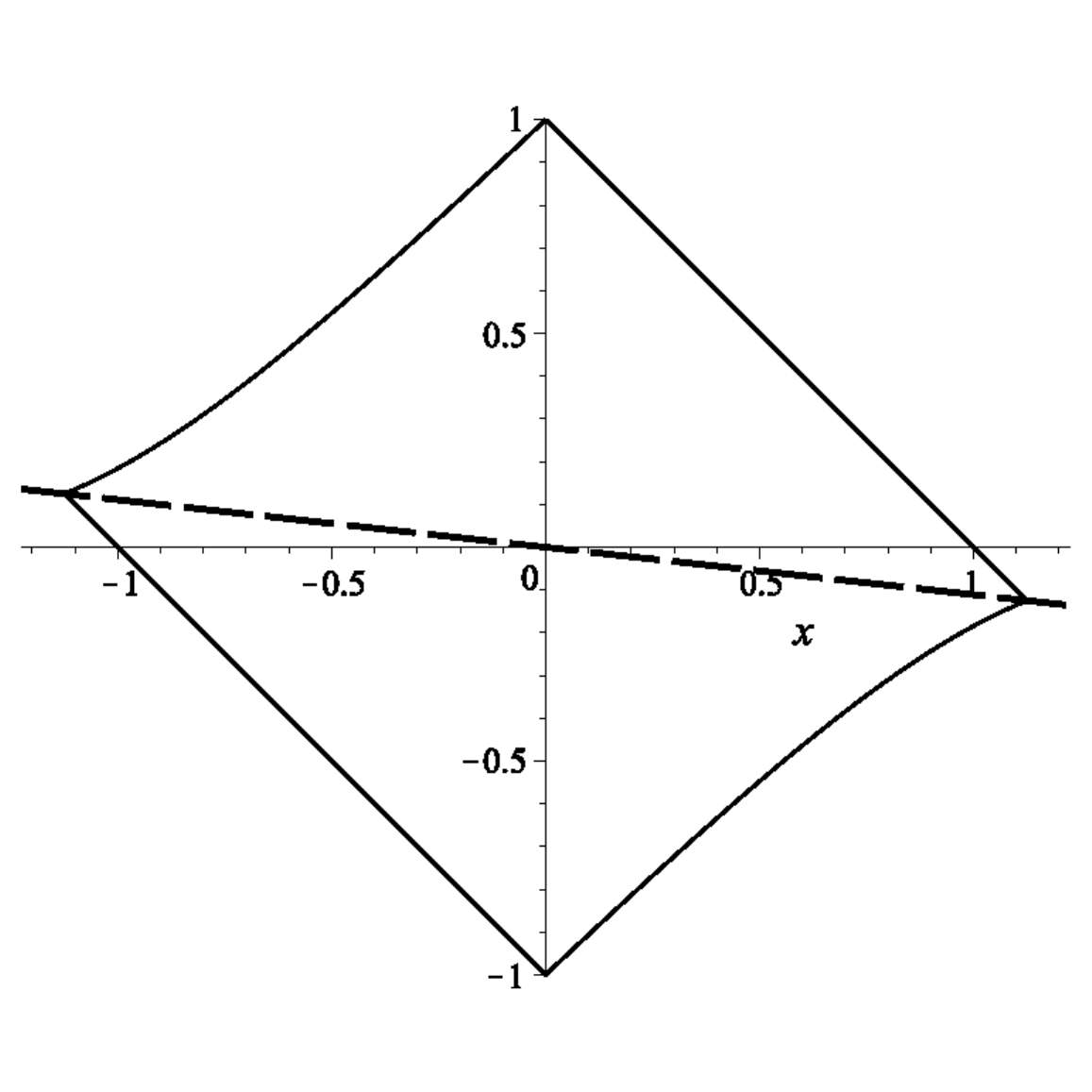}\\
\textit{ii)}
\end{center}
\end{minipage}
}
\caption{Stability domains of the trinomial $f(z)=z^n+az^{n-1}+b$ in the $(a,b)$-plane: 
{\it i)} even $n$, {\it ii)} odd $n$. The line $b=a/n$ is in bold, the line 
$b=-a/n$ is dashed.}
\label{fig:intersections}
\end{figure}

\begin{theorem}
All zeros of the quadrinomial $p(z)=1+\kappa(z+z^{N-1})+z^N$ lie on the unit circle if and
only if there hold the inequalities 
\[
-1\le\kappa\le
\left\{\begin{array}{ll}
     1 & \text{ if $N$ is {\it even}}, \\
     N/(N-2) & \text{ if $N$ is {\it odd}}.
\end{array}\right.
\]
\end{theorem}
\begin{proof}
Consider the quadrinomial $p(z)=1+\kappa(z+z^{N-1})+z^N$ and its derivative 
\[
p'(z)=N\left(\frac{\kappa}{N}+\frac{\kappa}{N}(N-1)z^{N-2}+z^{N-1}\right).
\]
Let us determine when all zeros of the trinomial $p'(z)$ belong to $\overline{\mathbb{D}}$. 
We apply Theorem~\ref{th:intro-2} with $n=N-1$, $b=\kappa/(n+1)$, $a=\kappa n/(n+1)$,
noting that $b=a/n$ and $a+b=\kappa$. The stability conditions for odd $n$ (i.e. N is even) can be 
written as $|a+b|\le1$ or $|\kappa|\le1$. 

Let $n$ be even. Then the intersection points of the line $b=a/n$ with the boundary of the 
stability domain are defined by the equations $nb+b=-1$ and $b=nb-1$, whence
\[
-\frac{1}{n+1}\le b\le\frac{1}{n-1}
\quad\text{or}\quad
-1\le\kappa\le\frac{n+1}{n-1}
\]
(see Figure~\ref{fig:intersections}.\textit{i)}). Returning from $n$ to $N$, we arrive at the
conclusion of the theorem. 
\end{proof}
\begin{theorem}
All zeros of the quadrinomial $q(z)=1+\kappa(z-z^{N-1})-z^N$ lie on the unit circle if and
only if there hold the inequalities 
\[
-N/(N-2)\le\kappa\le
\left\{\begin{array}{ll}
     1 & \text{ if $N$ is {\it odd}}, \\
     N/(N-2) & \text{ if $N$ is {\it even}}. 
\end{array}\right.
\]
\end{theorem}
\begin{proof}[Proof\nopunct]
goes as in the previous theorem. Figure~\ref{fig:intersections}.\textit{ii)}
is the illustration to the proof. 
\end{proof}
Next, we obtain factorization formulas for the quadrinomials 
$p(z)=1+\kappa(z+z^{N-1})+z^N$ at the limiting values of the parameter $\kappa$. The
well-known formulas \cite[1.396]{GR07} are 
\[
1-z^n=\left\{\begin{array}{ll}
    \displaystyle
    (1-z)\prod_{j=1}^{(n-1)/2}\left[1+z^2-2z\cos\frac{2\pi j}{n}\right] & 
        \text{ if $n$ is {\it odd}}, \\
    \displaystyle
    (1-z^2)\prod_{j=1}^{(n-2)/2}\left[1+z^2-2z\cos\frac{2\pi j}{n}\right] & 
        \text{ if $n$ is {\it even}},
\end{array}\right.
\]
\[
1+z^n=(1+z)\prod_{j=1}^{(n-1)/2}\left[1+z^2+2z\cos\frac{2\pi j}{n}\right]
\quad\text{ if $n$ is {\it odd}}.
\]
Using these identities, we easily derive factorization formulas for the polynomials $p(z)$
when $\kappa=-1$; $\kappa=1$ and $N$ is {\it even}. Let us write them in terms of the zeros
of the Chebyshev polynomial $U_{N-2}(x)$, while taking into account that 
$\cos2\alpha=2\cos^2\alpha-1$. If $\{\mu_j\}_{j=1}^{[(N-3)/2]}$ is the set of positive roots
of the equation $U_{N-2}(x)=0$ and $\beta_j=1-2(\mu_j)^2$, then we obtain the following
factorizations:
\begin{enumerate}
    \item[i)] $\kappa=-1$, $N$ is {\it odd},
    \[
    p(z)=(1+z)(1-z)^2\prod_{j=1}^{(N-3)/2}[1+z^2+2z\beta_j],
    \]
    \item[ii)] $\kappa=-1$, $N$ is {\it even},
    \[
    p(z)=(1-z)^2\prod_{j=1}^{(N-2)/2}[1+z^2+2z\beta_j],
    \]
    \item[iii)] $\kappa=1$, $N$ is {\it even},
    \[
    p(z)=(1+z)^2\prod_{j=1}^{(N-2)/2}[1+z^2-2z\beta_j].
    \]
\end{enumerate}

\begin{figure}[ht]
\includegraphics[width=0.5\textwidth]{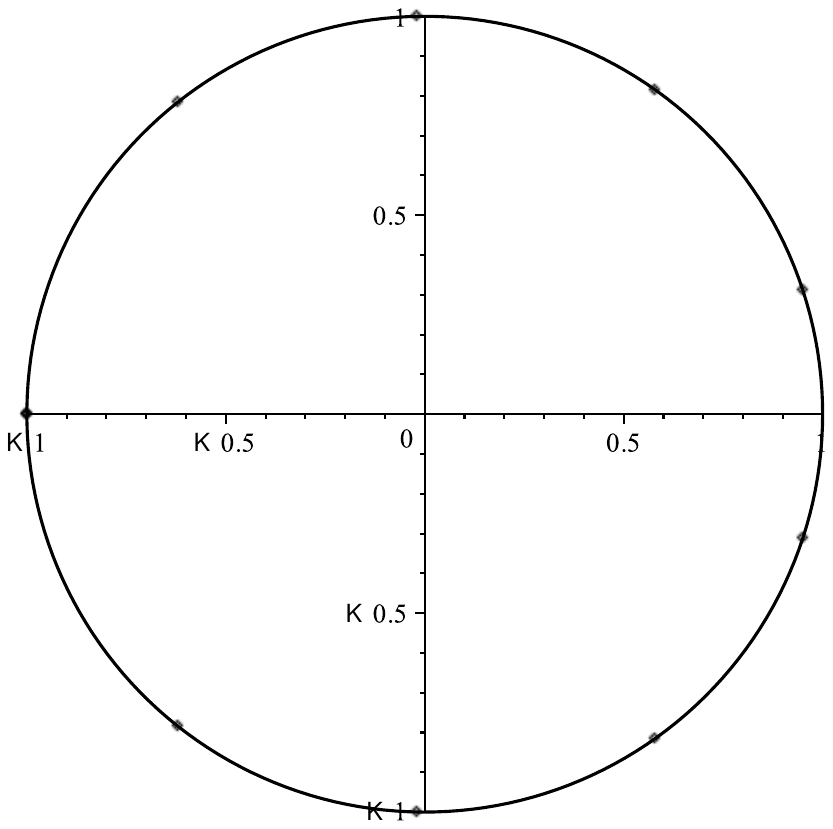}
\vspace{-3cm}
\caption{Location of the points $(-1,0)$ and $(\gamma_j,\pm\sqrt{1-(\gamma_j)^2})$,
$j=1,\ldots,(N-3)/2=4$ for $N=11$, on the unit circle.}
\label{fig:points}
\end{figure}

The factorization formula for the case $\kappa=N/(N-2)$, $N$ is {\it odd}, stands apart: 
instead of the zeros of the Chebyshev polynomial $U_{N-2}(x)$, it contains zeros of the 
polynomial's derivative $U'_{N-2}(x)$.
\begin{theorem}\label{th:p-factor}
Let $\{\nu_j\}_{j=1}^{(N-3)/2}$ be the set of positive roots of the equation $U'_{N-2}(x)=0$,
$\gamma_j=1-2(\nu_j)^2$, 
\[
p(z)=1+\frac{N}{N-2}(z+z^{N-1})+z^N.
\]
Then
\[
p(z)=(1+z)^3\prod_{j=1}^{(N-3)/2}[1+z^2-2z\gamma_j].
\]
\end{theorem}
\begin{proof}
Denote 
\[
\hat{p}(z)=(1+z)^3\prod_{j=1}^{(N-3)/2}[1+z^2-2z\gamma_j]
\]
and compute
\begin{align*}
p(e^{it})&=e^{i\frac{N}{2}t}\left(e^{-i\frac{N}{2}t}+e^{i\frac{N}{2}t}\right)
+\frac{N}{N-2}e^{it}e^{i\frac{N-2}{2}t}\left(e^{-i\frac{N-2}{2}t}+
e^{i\frac{N-2}{2}t}\right) \\
&=\frac{2}{N-2}e^{i\frac{N}{2}t}\left((N-2)\cos\frac{N}{2}t+N\cos\frac{N-2}{2}t\right); \\
\hat{p}(e^{it})&=2^{\frac{N+3}{2}}e^{i\frac{N}{2}t}\cos^3\frac{t}{2}
\prod_{j=1}^{(N-3)/2}[\cos t-\gamma_j].
\end{align*}
If $p(e^{it})=\hat{p}(e^{it})$ holds for every $t\in\mathbb{R}$, then for all $z$ we have
$p(z)=\hat{p}(z)$. Making the change of variables $\vartheta=(t-\pi)/2$, we rewrite
$p(e^{it})=\hat{p}(e^{it})$ equivalently as 
\[
\frac{1}{2\sin^3\vartheta}(N\sin(N-2)\vartheta-(N-2)\sin N\vartheta)
=2^{N-2}(N-2)\prod_{j=1}^{(N-3)/2}[\cos^2\vartheta-(\nu_j)^2].
\]
Since $U_{N-2}(x)=2^{N-2}x^{N-2}-\ldots$, $U'_{N-2}(x)=2^{N-2}(N-2)x^{N-3}-\ldots$, and
the function $U'_{N-2}(x)$ is even when $N$ is {\it odd}, there is the representation
\[
U'_{N-2}(\cos\vartheta)=2^{N-2}(N-2)\prod_{j=1}^{(N-3)/2}[\cos^2\vartheta-(\nu_j)^2].
\]
Next we use the formula \cite[Lemma~2]{DSS22}
\[
U'_k(x)=\frac{1}{2(1-x^2)}((k+2)U_{k-1}(x)-kU_{k+1}(x)),
\]
from which it follows that 
\[
U'_{N-2}(\cos\vartheta)=\frac{1}{2\sin^3\vartheta}(N\sin(N-2)\vartheta-(N-2)\sin N\vartheta).
\]
This gives the conclusion of the theorem.
\end{proof}
The factorization formulas for the quadrinomial $q(z)=1+\kappa(z-z^{N-1})-z^N$ are obtained
similarly. 
\begin{theorem}
Let $\{\nu_j\}_{j=1}^{(N-3)/2}$ be the set of positive roots of the equation $U'_{N-2}(x)=0$,
$\gamma_j=1-2(\nu_j)^2$, $q(z)=1+\kappa(z-z^{N-1})-z^N$. Then
\begin{enumerate}
    \item[i)] $\kappa=-N/(N-2)$, $N$ is {\it odd},
    \[
    q(z)=(1-z)^3\prod_{j=1}^{(N-3)/2}[1+z^2+2z\gamma_j],
    \]
    \item[ii)] $\kappa=-N/(N-2)$, $N$ is {\it even},
    \begin{equation}\label{eq:q-K_neg-N_even}
    q(z)=(1+z)(1-z)^3\prod_{j=1}^{(N-4)/2}[1+z^2+2z\gamma_j],
    \end{equation}
    \item[iii)] $\kappa=N/(N-2)$, $N$ is {\it even},
    \[
    q(z)=(1-z)(1+z)^3\prod_{j=1}^{(N-4)/2}[1+z^2-2z\gamma_j].
    \]
\end{enumerate}
\end{theorem}
\begin{remark}
In the case when $\kappa=1$ and $N$ is {\it odd}, the factorization is
\[
q(z)=(1-z)(1+z)^2\prod_{j=1}^{(N-3)/2}[1+z^2+2z\beta_j],
\]
where $\beta_j=1-2(\mu_j)^2$, $\{\mu_j\}_{j=1}^{(N-3)/2}$ is the set of positive roots of 
the equation $U_{N-2}(x)=0$.
\end{remark}
When $\kappa=\pm1$, the zeros of the polynomials $p(z)$ and $q(z)$ are evenly spaced around 
the circle, that is, they form a set of the roots of unity. The situation is different in the
case $\kappa=N/(N-2)$, as illustrated in Figure~\ref{fig:points}.

The sequence $\{\arccos\gamma_j\}_{j=1}^{(N-3)/2}$ is {\it not} arithmetic. For $N=11$,
we have
\begin{center}
\begin{tabular}{ l l }
    $\arccos\gamma_1=0.3173\ldots$, & \\ 
    $\arccos\gamma_2=0.9527\ldots$, & 
        $\arccos\gamma_2-\arccos\gamma_1=0.6354\ldots\ne2\arccos\gamma_1$; \\
    $\arccos\gamma_3=1.5911\ldots$, & $\arccos\gamma_3-\arccos\gamma_2=0.6383\ldots$; \\
    $\arccos\gamma_4=2.2398\ldots$, & $\arccos\gamma_4-\arccos\gamma_3=0.6487\ldots$.
\end{tabular}
\end{center}
\section{Factorization of the derivative of the Fej\'er polynomial}
Let $f(z)=\sum_{j=1}^n z^j$, then
\[
\frac{f(ze^{i\alpha_k})-f(ze^{-i\alpha_k})}{z(e^{i\alpha_k}-e^{-i\alpha_k})}
=\frac{1-(-1)^k z^{n+1}}{z^2+1-2z\cos\alpha_k},
\]
where $\alpha_k=k\pi/(n+1)$, $k=1,\ldots,n$; that is, $f(z)\in\mathcal{P}_n$. It follows from
Theorem~\ref{th:Suffridge} that the Fej\'er polynomial 
$\sigma_n(z)=\sum_{j=1}^n\left(1-\frac{j-1}{n}\right)z^j$ is univalent in $\mathbb{D}$.

Apparently the first proof of the univalence of the Fej\'er polynomial was given 
in~\cite{EE37}. It is also shown there that 
\[
\sigma'_N(z)=\frac1{(1-z)^3}\left(1-\frac{N+2}{N}z+\frac{N+2}{N}z^{N+1}-z^{N+2}\right).
\]
Therefore, for odd $N$ there holds
\[
(1+z)^3\sigma'_{N-2}(-z)=1+\frac{N}{N-2}(z+z^{N-1})+z^N,
\]
and hence Theorem~\ref{th:p-factor} implies that when $N$ is {\it odd},
\begin{equation}\label{eq:sigma-odd_N}
\sigma'_N(z)=\prod_{j=1}^{(N-1)/2}[1+z^2+2z\tilde{\gamma}_j],
\end{equation}
where $\tilde{\gamma}_j=1-2(\tilde{\nu}_j)^2$, $\{\tilde{\nu}_j\}_{j=1}^{(N-1)/2}$ is the
set of positive roots of the equation $U'_N(x)=0$. That is, preimages of cusps in the image
of the unit circle under the mapping by the Fej\'er polynomial form the set 
$\{e^{\pm i\tilde{\gamma}_j}\}_{j=1}^{(N-1)/2}$. Using formula~\eqref{eq:q-K_neg-N_even}, 
we easily obtain the following analogue of formula~\eqref{eq:sigma-odd_N} for 
{\it even} $N$:
\begin{equation}\label{eq:sigma_even_N}
\sigma'_N(z)=(1+z)\prod_{j=1}^{(N-2)/2}[1+z^2+2z\tilde{\gamma}_j].
\end{equation}
\begin{remark}
Consider {\it J.~Alexander}'s polynomial~\cite{JA15} $w_N(z)=\sum_{j=1}^N\frac1j z^j$. Denote
by $\{\tilde{\mu}_j\}_{j=1}^{[(N-1)/2]}$ the set of positive roots of the equation 
$U_{N-1}(x)=0$, $\tilde{\beta}_j=1-2(\tilde{\mu}_j)^2$. Then the formulas 
\[
w'_N(z)=\left\{\begin{array}{ll}
    \displaystyle
    \prod_{j=1}^{(N-1)/2} [1+z^2+2\tilde{\beta}_j z] & \text{if $N$ is {\it odd}}, \\
    \displaystyle
    (1+z)\prod_{j=1}^{(N-2)/2} [1+z^2+2\tilde{\beta}_j z] & \text{if $N$ is {\it even}}
\end{array}\right.
\]
hold true and may be considered as the analogue of formulas~\eqref{eq:sigma-odd_N} and
\eqref{eq:sigma_even_N}. 
\end{remark}
\section{Quadrinomial $p(z)=1+\frac{N}{N-2}(z+z^{N-1})+z^N$ and some univalent polynomials}
Let again $N$ be {\it odd}. Consider the polynomial 
\[
\tilde{p}(z)=\frac{z}{(1+z)^2}\left(1+\frac{N}{N-2}(z+z^{N-1})+z^N\right).
\]
All zeros of this polynomial are simple and located on the unit circle (except for the
trivial root). This polynomial can be represented by its coefficients with the use of the 
summation formula for a geometric sequence and its derivative, namely
\[
\tilde{p}(z)=\sum_{j=1}^{(N-1)/2}(-1)^{j-1}\left(1-2\frac{j-1}{N-2}\right)(z^j+z^{N-j}).
\]
Then, let us apply the Suffridge transformation 
\[
F(z)=\frac{N}{N-1}\tilde{p}(z)-\frac{1}{N-1}z\tilde{p}'(z),
\]
or
\begin{equation}\label{eq:F-univalent}
F(z)=\sum_{j=1}^{(N-1)/2}(-1)^{j-1}\left(1-2\frac{j-1}{N-2}\right)
\left(\frac{N-j}{N-1}z^j+\frac{j}{N-1}z^{N-j}\right).
\end{equation}
\begin{theorem}\label{th:univalence}
Polynomial~\eqref{eq:F-univalent} is univalent in $\mathbb{D}$.
\end{theorem}
\begin{proof}
Let $\alpha_k=k\pi/N$, $k=1,\ldots,N$. Denote
\[
\varphi_k(z)=\frac{\tilde{p}(ze^{i\alpha_k})-\tilde{p}(ze^{-i\alpha_k})}
{z(e^{i\alpha_k}-e^{-i\alpha_k})}.
\]
A straightforward computation yields the relation
\[
\varphi_k(z)=(1-(-1)^k z^{N+2})+\frac{2N}{N-2}\cos{\alpha_k}(z-(-1)^k z^{N+1})+
\frac{N+2}{N-2}(z^2-(-1)^k z^N).
\]
If the polynomials $\varphi_k(z)$, $k=1,\ldots,N$, do not have zeros in $\mathbb{D}$, then
by Theorem~\ref{th:Suffridge} polynomial~\eqref{eq:F-univalent} is univalent. We will prove
a more precise result, namely, that {\it all zeros of the polynomials $\varphi_k(z)$
are located on the unit circle}. 

Compute
\begin{align*}
\varphi_k&(e^{it})=\frac{2}{N-2}e^{i\frac{N+2}{2}t} \\
&\times\left\{\begin{array}{l}
    \displaystyle
    (N-2)\cos\frac{N+2}{2}t+2N\cos\alpha_k\cos\frac{N}{2}t+(N+2)\cos\frac{N-2}{2}t \\ 
    \hfill\text{if $k$ is {\it odd}}, \\
    \displaystyle
    -i\left((N-2)\sin\frac{N+2}{2}t+2N\cos\alpha_k\sin\frac{N}{2}t+(N+2)\sin\frac{N-2}{2}t\right)\qquad\\ 
    \hfill\text{if $k$ is {\it even}}.
\end{array}\right.
\end{align*}
We need to show that the equation $\varphi_k(e^{it})=0$ has precisely $N+2$ roots on the 
interval $[0,2\pi]$. Cases `$k$ is {\it odd}' and `$k$ is {\it even}' are treated analogously.
Let us deal with the case when $k$ is {\it odd}.

Denote
\[
C(\vartheta)=(N-2)\cos(N+2)\vartheta+2N\cos\alpha_k\cos N\vartheta+(N+2)\cos(N-2)\vartheta.
\]
Since $C(\pi/2)=0$ and $C(\vartheta)=-C(\pi-\vartheta)$, then we must show that the 
trigonometric polynomial $C(\vartheta)$ has precisely $(N+1)/2$ roots on the interval 
$(0,\pi/2)$. Let $\vartheta_j=j\pi/N$, $j=0,\ldots,(N-1)/2$. Then 
$C(\vartheta_j)=(-1)^j\,2N(\cos2\vartheta_j+\cos\alpha_k)$. If $j_0=(N-k)/2$, then 
$C(\vartheta_{j_0})=0$. However, the root $\vartheta_{j_0}$ is multiple. Indeed,
\[
C'(\vartheta)
=-2(N^2-4)\sin N\vartheta\left(\cos2\vartheta+\frac{N^2}{N^2-4}\cos\alpha_k\right),
\]
that is, $C'(\vartheta_{j_0})=0$. Consider the sequence of signs 
\[
\left\{\textnormal{sign}\,C(\vartheta_j)\right\},
\quad j\in\left\{0,\ldots,\frac{N-1}{2}\right\}\setminus\{j_0\}.
\]
The signs in this sequence---consisting of $(N-1)/2$ elements---strictly alternate: 
$\left\{+,-,+,\ldots,(-1)^\frac{N+1}{2}\right\}$. This means that the trigonometric polynomial
$C(\vartheta)$ has $(N-3)/2$ simple roots on $(0,\pi/2)$. We add the multiple root to them
and obtain precisely $(N+1)/2$ roots of the polynomial $C(\vartheta)$ on the interval
$(0,\pi/2)$, which proves the theorem.
\end{proof}
Thus, polynomial~\eqref{eq:F-univalent} is univalent. Three more univalent polynomials can be
constructed in a similar way. Denote $-F(-z)=F_N^{(1)}(z)$. Note that when $N$ is {\it even},
\[
\frac{z}{(1+z)^2}\left(1+\frac{N}{N-2}(z-z^{N-1})-z^N\right)
=\sum_{j=1}^{(N-2)/2}\left(1-2\frac{j-1}{N-2}\right)(z^j-z^{N-j}). 
\]
Proceeding next as in the proof of Theorem~\ref{th:univalence}, we finally obtain
\[
F_N^{(1)}(z)=\sum_{j=1}^{(N-1)/2}\left(1-2\frac{j-1}{N-2}\right)
\left(\frac{N-j}{N-1}z^j-\frac{j}{N-1}z^{N-j}\right),\quad\text{$N$ is {\it odd}},
\]
\[
F_N^{(2)}(z)=\sum_{j=1}^{(N-1)/2}\left(1-2\frac{j-1}{N-2}\right)
\left(\frac{N-j}{N-1}z^j+\frac{j}{N-1}z^{N-j}\right),\quad\text{$N$ is {\it odd}},
\]
\[
F_N^{(3)}(z)=\sum_{j=1}^{(N-2)/2}\left(1-2\frac{j-1}{N-2}\right)
\left(\frac{N-j}{N-1}z^j-\frac{j}{N-1}z^{N-j}\right),\quad\text{$N$ is {\it even}},
\]
\[
F_N^{(4)}(z)=\sum_{j=1}^{(N-2)/2}\left(1-2\frac{j-1}{N-2}\right)
\left(\frac{N-j}{N-1}z^j+\frac{j}{N-1}z^{N-j}\right),\quad\text{$N$ is {\it even}}.
\]
\begin{figure}[ht]
{\footnotesize
\begin{minipage}{.35\textwidth}
\begin{center}
\includegraphics[width=\textwidth]{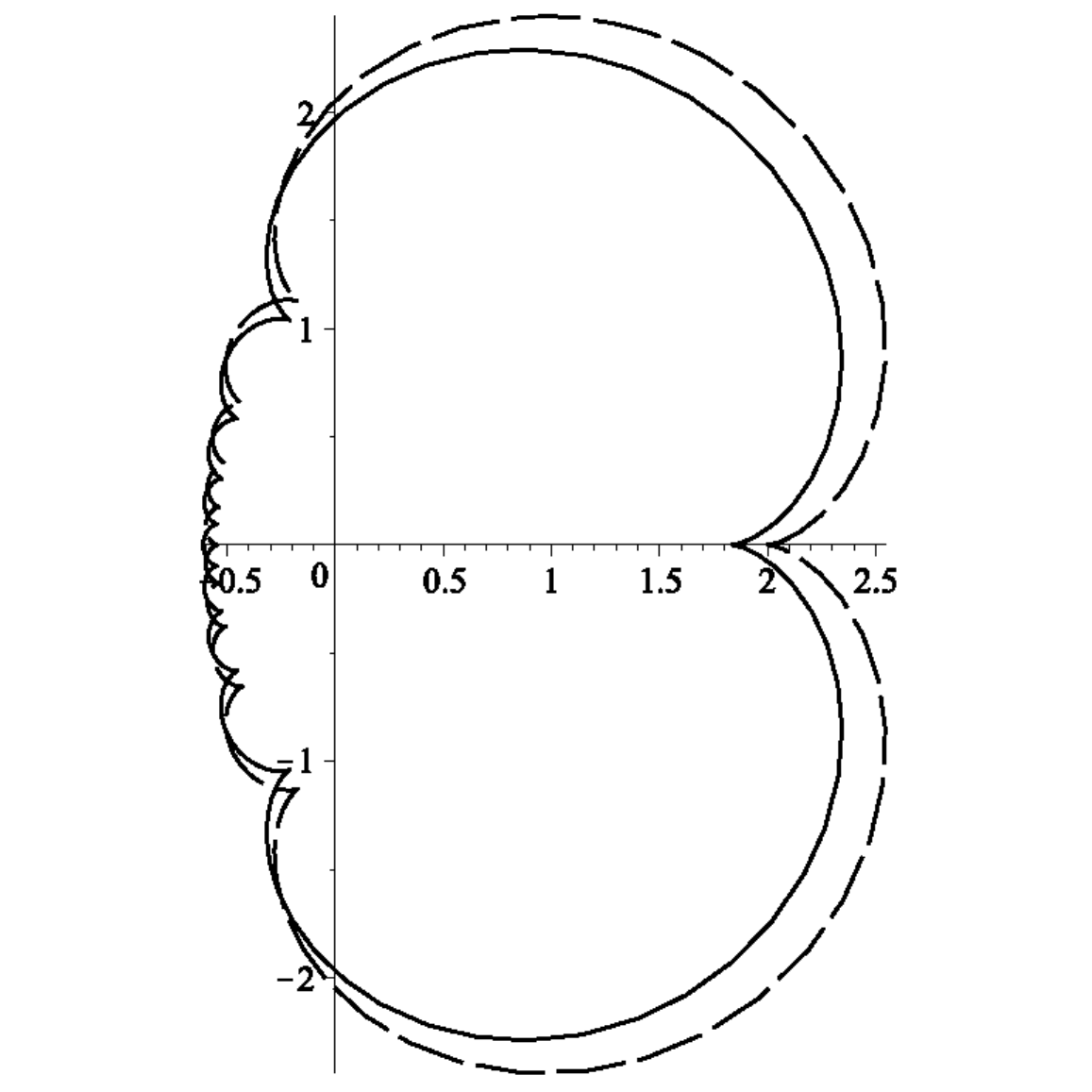}\\
\textit{i)}
\end{center}
\end{minipage}
\hspace{12pt}
\begin{minipage}{.35\textwidth}
\begin{center}
\includegraphics[width=\textwidth]{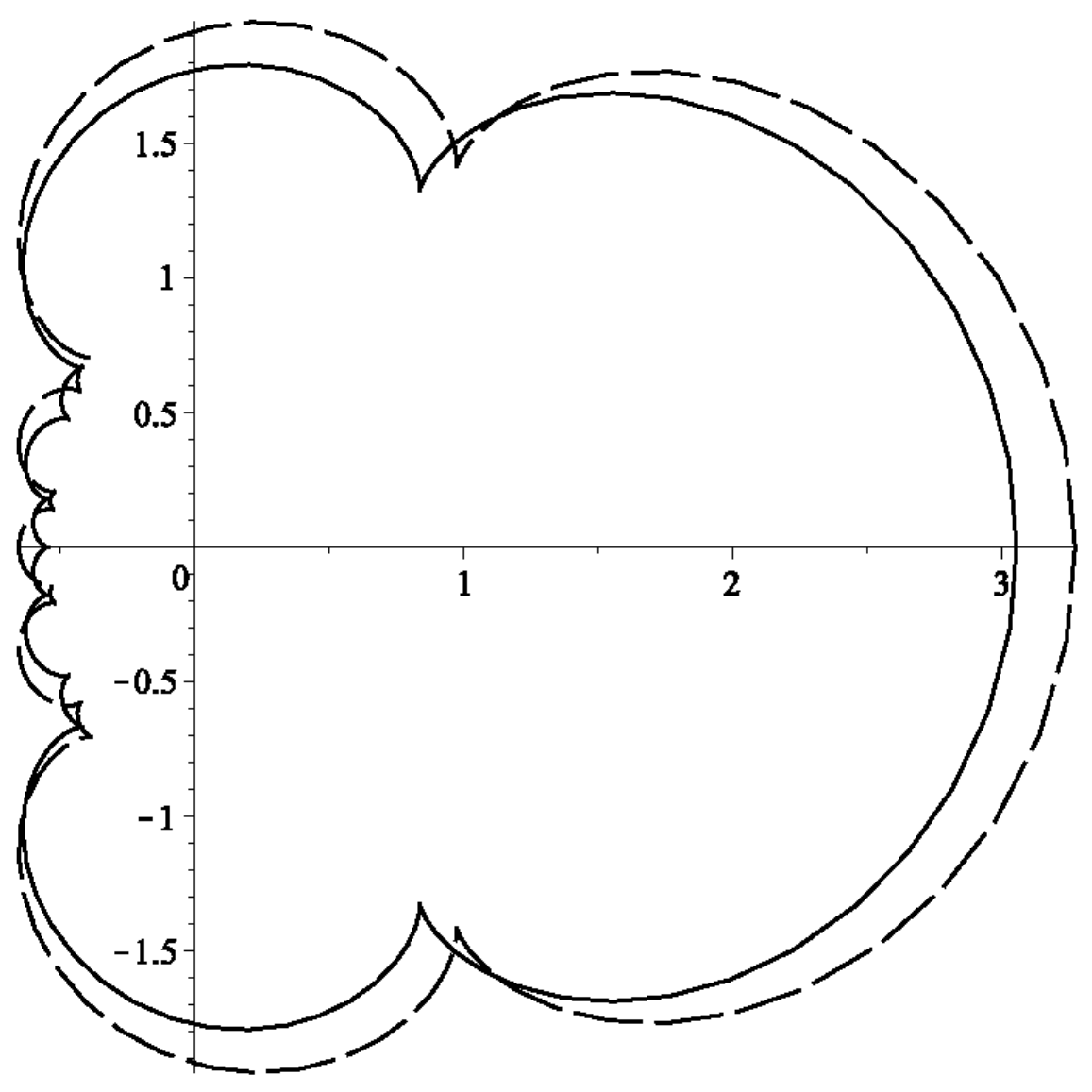}\\
\textit{ii)}
\end{center}
\end{minipage}
}
\caption{Graphs of the images of the unit circle: {\it i)} under the mappings $F_{11}^{(1)}(z)$
and $F_{12}^{(3)}(z)$ (dashed); {\it ii)} under the mappings $F_{11}^{(2)}(z)$ and
$F_{12}^{(4)}(z)$ (dashed).}
\end{figure}
Among the univalent in $\mathbb{D}$ polynomials, the class of the polynomials such that 
all zeros of their derivative lie on the unit circle is of great importance. 
{\it St.~Ruscheweyh} called these polynomials quasi-extremal~\cite[p.~284]{GS02}. The
polynomials $F_N^{(s)}(z)$ built above, $s=1,2,3,4$, are quasi-extremal. Let us demonstrate
this on the example of polynomial~\eqref{eq:F-univalent}.
\begin{theorem}
All zeros of the derivative of polynomial~\eqref{eq:F-univalent} lie on the unit circle.
\end{theorem}
\begin{proof}
Let 
\[
\tilde{p}(z)=\sum_{j=1}^{(N-1)/2}(-1)^{j-1}\left(1-2\frac{j-1}{N-2}\right)(z^j+z^{N-j})
\]
and
\[
F(z)=\frac{N}{N-1}\tilde{p}(z)-\frac{1}{N-1}z\tilde{p}'(z),
\]
$N$ is {\it odd}. Consider the polynomial $F'(z)=\tilde{p}'(z)-\frac{1}{N-1}z\tilde{p}''(z)$.
After the transformations we arrive at
\[
F'(z)=\frac{1}{(N-1)(N-2)(1+z)^4}W(z),
\]
where $W(z)=(N-1)(N-2)(1+z^{N+2})+2(N-2)(N+2)(z+z^{N+1})+(N+1)(N+2)(z^2+z^N)$. It is not 
difficult to check that the polynomial $W(z)$ has the root $z_0=-1$ of multiplicity five.
Let us show that the other $N-3$ roots lie on the unit circle too.

Calculate
\begin{align*}
W(e^{it})=2e^{i\frac{N+2}{2}t}&\left[(N-1)(N-2)\cos\frac{N+2}{2}t\right. \\
&\qquad\left.+2(N-2)(N+2)\cos\frac{N}{2}t+(N+1)(N+2)\cos\frac{N-2}{2}t\right].
\end{align*}
Denote 
\begin{align*}
C(\vartheta)=(N-1)&(N-2)\cos(N+2)\vartheta \\
&+2(N-2)(N+2)\cos N\vartheta+(N+1)(N+2)\cos(N-2)\vartheta. 
\end{align*}
Since $C(\pi/2)=0$ and $C(\vartheta)=-C(\pi-\vartheta)$, we
need to show that the trigonometric polynomial $C(\vartheta)$ has precisely $(N-3)/2$ roots
on the interval $(0,\pi/2)$. Let $\vartheta_j=j\pi/N$, $j=0,\ldots,(N-1)/2$. Then
$C(\vartheta_j)=(-1)^j((N^2+2)\cos2\vartheta_j+2)=(-1)^j\,4N^2\cos^2\vartheta_j$. Thus the
function $C(\vartheta)$ alternates its sign $(N-3)/2$ times on the interval 
$\left[\frac{\pi}{N},\frac{N-1}{2}\frac{\pi}{N}\right]$. Hence this function has precisely
$(N-3)/2$ roots on this interval. The theorem is proved.
\end{proof}

\bibliographystyle{unsrt}
\bibliography{main.bib}
\end{document}